\input amstex
\documentstyle{amsppt}
\voffset=-3pc
\def\bC{\Bbb C}
\def\bB{\Cal B}

\def\bU{\Cal U}
\def\bK{\Bbb K}
\def\os{\overline S}

\def\oU{\overline {\Cal U}}

\def\ep{\epsilon}
\def\Rsa{R_{\text{sa}}}
\def\asa{A_{\text{sa}}}
\def\csa{C_{\text{sa}}}
\def\tA{\widetilde A}
\def\Asa{\widetilde A_{\text{sa}}}
\def\Csa{\widetilde C_{\text{sa}}}
\magnification=\magstep1
\parskip=6pt
\NoBlackBoxes
\topmatter
\title Large C*--algebras of Universally Measurable Operators
\endtitle  
\author Lawrence G.~Brown
\endauthor
%\endtopmatter

\abstract{For a $C^*$--algebra $A$, G. Pedersen defind the concept of universal measurability for self--adjoint elements of $A^{**}$, the enveloping von Neumann algebra of $A$. Although he was unable to show that $\bU$, the set of universally measurable elements, is a Jordan algebra, he showed that $\bU$ contains a large Jordan algebra, $\bU_0$. We show by means of a $2\times 2$ matrix trick that $\bU_0$ is in fact the real part of a $C^*$--algebra. If it is ever shown that $\bU$ is always a Jordan algebra, then the same trick will show that $\bU$ is the real part of a $C^*$--algebra.}
\endabstract
\endtopmatter

Pedersen used the theory of semicontinuous elements of $A^{**}$ to define universal measurability. The basic reference for semicontinuity theory is \cite{AP}. An element of $A^{**}$  is {\it strongly lower semicontinuous} if it is in $((A_{\text{sa}})^m)^-$.  Here, for a subset $R$ of $A^{**}$, $\Rsa$ denotes $\{x\in R: x^*=x\}$, $(\Rsa)^m$ denotes the set of limits of bounded increasing nets from $\Rsa$, and $^-$ denotes norm closure. Also, the set of {\it middle lower semicontinuous} elements is $(\Asa)^m$ and the set of {\it weakly lower 
semicontinuous} elements is $((\Asa)^m)^-$, where $\tA=A+\bC 1$ and 1 is the identity of $A^{**}$. And $h$ is upper semicontinuous in any of the senses if and only if $-h$ is lower semicontinuous. Then for $x$ in $A^{**}_{\text sa}$, $x\in \bU=\bU(A)$ ([P1, p.435] or [P2, p.104]) if and only if for any $\ep>0$ and any state $\varphi$ of $A$, there are a strongly lower semicontinuous $h$ and a strongly upper semicontinuous $k$ such that
$$ k\leq x\leq h\quad\text{and }\quad \varphi(h-k)<\ep.$$
Also $x$ is in $\bU_b$ if and only if there is a constant $\gamma$, depending only on $x$, such that $h$ and $k$ above can always be taken with $\|h\|, \|k\|\leq \gamma$. Finally, $\bU_0$ ([P1, p.438]) is the monotone sequential closure of $\bU_b$, and $\overline S$ ([C, \S 2.2]) is the norm closed real vector space generated by $((\asa)^m)^-$. 

It was shown by F. Combes in [C, Prop. 2.2.15] that $\overline S$ is a Jordan algebra. (Actually [C, Prop. 2.2.15] is stated only for the unital case. However, Combes' general definition of semicontinuity in [C, \S 2.3] makes it clear that the result holds in general, provided it is shown that Combes' definition agrees with the Akemann--Pedersen concept of strong semicontinity. This last can readily be deduced from the results of \cite{AP}, or one can cite [B, Prop. 2.14].) Pedersen proved in [P1, Lemma 3.5 and Prop. 3.6] that $\overline \bU_b$ and $\bU_0$ are Jordan algebras and that $\bU_0\subset \bU$. It is easy to see that
$$\overline S\subset\overline \bU_b\subset \bU_0$$
and that all semicontinuous elements of any type, are contained in $\overline S$.

We will show that each of these Jordan algebras is the real part of a $C^*$--algebra. In other words $\overline S_{\bC}, \overline \bU_{b\bC}$, and $\bU_{0\bC}$ are $C^*$--algebras, where $R_\bC$ denotes $R+iR$ for any $R\subset A^{**}_{\text{sa}}$. It follows that $\bB$, the monotone sequentially closed $C^*$--algebra generated by the semicontinuous elements, is contained in $\bU_{\bC}$. In [C, Def. 2.4.1] Combes defined a Borel element to be a member of the monotone sequential closure of $\overline S$. Since $\os_{\bC}$ is a $C^*$--algebra and since Pedersen proved in [P2, Theorem 4.5.4] that the monotone sequential closure of a $C^*$--algebra is a $C^*$--algebra (an extension of an earlier Jordan algebra result of R. Kadison \cite{K}), we see that $\bB$ is just the complexification of Combes' set of Borel elements. It seems that $\bB$ is the most natural non--commutative analogue of the algebra of bounded Borel functions on a locally compact Hausdorff space. Note that (by the results of \cite{AP}) $\bB$ contains all open or closed projections and all (left, right, or two--sided) multipliers or quasi--multipliers of $A$.

Many applications of $\bU$ are based on the fact [P1, Theorem 3.8] that the atomic representation is isometric on $\bU$ (in particular, elements of $\bU$ are determined by their atomic parts) or the fact (told to the author by Pedersen) that elements of $\bU$ satisfy the barycenter formula. For example, this is why the authors of \cite{BW} needed to know, in connection with [BW, Lemma 3.7], that every self--adjoint element of the $C^*$--algebra generated by $A$ and a single closed projection is universally measurable. This could be established by {\it ad hoc} methods, but those methods wouldn't work, for example, for the $C^*$--algebra generated by $A$ and two closed projections. Since $\bB$ is a $C^*$--algebra, it seems to be large enough to cover future applications of this sort.

The $\bU$ part of the next proposition follows from [AAP, Lemma 2.1], and the proofs of the other parts are similar.

\proclaim{Proposition 1}Each of the spaces $\overline{S}_\bC$, $\bU_{b\bC}, \overline \bU_{b\bC}, \bU_{0\bC}$, and $ \bU_\bC$ is an $M(A)-M(A)-$ bimodule, where $M(A)$ is the multiplier algebra of $A$.
\endproclaim

\demo{Proof}For $T$ in $M(A)$, let $\theta_T(x)=T^*xT$. Since $\theta_T$ is positive, bounded, and carries $A$ into $A$, we see that $\theta_T$ preserves strong semicontinuity; and it is routine to show that $\theta_T$ maps each of the five spaces into itself. For example, if $h$ is in $\bU_b$, to approximate $T^*hT$ at a state $\varphi$, we approximate $h$ at the positive functional $\varphi(T^*\cdot T)$. If $\gamma$ is the constant appearing in the definition of $\bU_b$ for $h$, then the constant for $T^*hT$ will be $\gamma\|T\|^2$. Then for $h$ in one of the complexified spaces and $S, T$ in $M(A)$, it follows from polarization that $S^*hT$ is in the same space. Since $1\in M(A)$, the result follows.
\enddemo

Let $\{e_{ij}:i,j=1, \dots, n\}$ be the standard matrix units in $\Bbb M_n(A^{**})$. Thus for $x_{ij}$ in $A^{**}$, $(x_{ij})=\Sigma x_{ij}e_{ij}$. Note that $e_{ij}\in M(\Bbb M_n(A))$.

\proclaim{Lemma 2} For $x$ in $A^{**}_{\text sa}$, $x$ is in $\os(A),\bU_b(A), 
\oU_b(A), \bU_0(A)$, or $\bU(A)$ if and only if $xe_{11}$ is in $\os(\Bbb M_n(A)), \bU_b(\Bbb M_n(A)), \oU_b(\Bbb M_n(A)), \bU_0(\Bbb M_n(A))$, or $\bU(\Bbb M_n(A))$, respectively.
\endproclaim

\demo{Proof}Now we use two positive maps. Let $\theta(x)=xe_{11}$ and $\psi((a_{ij}))=a_{11}$. It is routine to see that each map preserves strong semicontinuity, and then to see that it preserves each of the five spaces. For example, to approximate $\psi(h)$, $h\in \bU(\Bbb M_n(A))$, at a state $\varphi$, we approximate $h$ at the state $\varphi\circ\psi$.
\enddemo

\proclaim{Proposition 3}If $x=(x_{ij})$ is in $\Bbb M_n(A^{**})$, then $x\in\os_\bC(\Bbb M_n(A))$, $\bU_{b\bC}(\Bbb M_n(A))$,
$\oU_{b\bC}(\Bbb M_n(A)), \bU_{0\bC}(\Bbb M_n(A))$, or $\bU_\bC(\Bbb M_n(A))$ if and only if each $x_{ij}$ is in $\os_\bC(A), \bU_{b\bC}(A)$,
 $\oU_{b\bC}(A)$, $\bU_{0\bC}(A)$, or $\bU_\bC(A)$, respectively.
\endproclaim

\demo{Proof}Combine the previous two results. Note that $x_{ij}e_{11}=e_{1i}xe_{j1}$ and $x=\Sigma e_{i1}(x_{ij}e_{11})e_{1j}$.
\enddemo

\proclaim{Theorem 4}Each of the spaces $\os_\bC, \oU_{b\bC}$, and $\bU_{0\bC}$ is a $C^*$--algebra. Moreover, $\bU_\bC$ is a $\bU_{0\bC}-\bU_{0\bC}$ -- bimodule.
\endproclaim

\demo{Proof}If $x$ is in $\bU_{0\bC}$, for example, then $y=\pmatrix 0 & x^*\\
x & 0\endpmatrix$ is in $\bU_{0\bC}(\Bbb M_2(A))$. Since $\bU_{0\bC}(\Bbb M_2(A))$ is a Jordan algebra, then $y^2=\pmatrix x^*x & 0\\
0 & xx^*\endpmatrix$ is in $\bU_{0\bC}(\Bbb M_2(A))$, and $x^*x\in \bU_{0\bC}$. By polarization, $x^*y\in \bU_{0\bC}$ wherever $x,y\in \bU_{0\bC}$.

For the second statement, note that $a^*ha\in \bU_0\subset \bU$, wherever $a\in \bU_{0\bC}$ and $h$ is strongly semicontinuous. It follows, as in the proof of Proposition 1, that if $x$ is in $\bU$, $\varphi$ is a state, and $\ep>0$, we can find $h, k$ in $\bU$ such that
$$k\leq a^*xa\leq h\quad\text{and}\quad \varphi(h-k)<\ep.$$
Since $h$ can be approximated from above and $k$ from below by appropriate semicontinuous elements, it follows that $a^*xa\in \bU$, Now the result follows from polarization, as above.
\enddemo

\proclaim{Corollary 5}If $\bB$ is the monotone sequentially closed $C^*$--algebra generated by $((\asa)^m)^-$, then $\bB \subset \bU_\bC$.
\endproclaim

Of course $\bB$ is also contained in $\bU_{0\bC}$.

\example{Remark 6}There is a canonical maximal $C^*$--algebra in $\bU_\bC$. Let $\bU_1=\{x\in A^{**}: x \bU_\bC+ \bU_\bC x\subset \bU_\bC\}$.
It is clear {\it a priori} that $\bU_1$ is a $C^*$--algebra and $\bU_\bC$ is a $\bU_1-\bU_1$-- bimodule. Since $\bU$ is strongly sequentially closed by [P1, Theorem 3.7], $\bU_1$ is double--strongly sequentially closed. By Theorem 4, $\bU_1$ contains $\bU_0$. If $C$ is any $C^*$--algebra in $\bU_\bC$ such that $((\asa)^m)^-\subset C$, then the proof of the second statement of Theorem 4 shows that $C\subset \bU_1$. The main purpose of this remark is to point out that there is no cost to using the bimodule approach in the effort to construct large $C^*$--algebras in $\bU_\bC$. 
\endexample

We think the next result could have been proved earlier from the barycenter formula and Choquet theory, but the following proof seems easier.

\proclaim{Proposition 7}The atomic representation is completely isometric on $\bU_\bC$. 
\endproclaim

\demo{Proof} If $x\in \Bbb M_n(\bU_\bC)$, consider $y=\pmatrix 0 & x\\
x^* & 0\endpmatrix$ in $\Bbb M_{2n}(\bU_\bC)$. Since $y\in \bU(\Bbb M_{2n}(A))$ by Proposition 3, Pedersen's result, [P1, Theorem 3.8], implies 
$$\|x\|=\|y\|=\|y_{\text{atomic}}\|=\|x_{\text{atomic}}\|.$$
\enddemo

It may be of interest to consider more general tensor products in relation to universal measurability. Note that if $C=A\otimes B$, an arbitrary $C^*$--tensor product, there are normal embeddings of $A^{**}$ and $B^{**}$ in $C^{**}$.  By 
the proof of Theorem 2 of [A], this leads  
to an injection of $A^{**}\underset{\text{alg}}\to\otimes B^{**}$ into $C^{**}$. The following result generalizes Proposition 3.

\proclaim{Proposition 8}If $C=A\otimes B$, a $C^*$--tensor product, then $\os_\bC (A)\underset{\text{alg}}\to\otimes\os_\bC (B)\subset \os_\bC (C)$, $ \bU_{b\bC}(A)\underset{\text{alg}}\to\otimes \bU_{b\bC}(B)\subset \bU_{b\bC}(C)$, $\oU_{b\bC}(A)\underset{\text{alg}}\to\otimes\oU_{b\bC}(B)\subset\oU_{b\bC}(C)$, $\bU_{0\bC}(A)\underset{\text{alg}}\to\otimes \bU_{0\bC}(B)\subset \bU_{0\bC}(C)$, and $\bU_\bC(A)\underset{\text{alg}}\to\otimes \bU_\bC(B)\subset \bU_\bC(C)$. Moreover, every slice map carries each of the five spaces for $C$ into the corresponding space for $A$ or $B$.
\endproclaim

\demo{Proof}First we show that $\os(A)\otimes 1\subset\os(C)$, etc. For $h\in((\asa)^m)^-$, it is easy to see that $h\otimes 1\in ((\csa)^m)^-$ if $h\geq 0$ and $h\otimes 1\in(\Csa)^m$ in general. This is already enough to yield the result for $\os$ and $\bU$. If $x\in \bU_b(A)$, in order to approximate $x\otimes 1$ from above, say, we first approximate $x+\|x\|\cdot 1$ from above by $h$. Then $x\otimes 1$ is approximated from above by $h\otimes 1 -\|x\|e_i$, where $e_i$ is taken from an approximate identity of $C$. If $\gamma$ is the constant for $x$ appearing in the definition of $\bU_b$, then the constant for $x\otimes 1$ is no more then $\gamma +2\|x\|$. The $\oU_b$-- and $\bU_0$--cases follow immediately.

To deal with $x\otimes y$, we write $x=x_0+\lambda\cdot 1$ and $y=y_0+\mu\cdot 1$ where $x_0, y_0\geq 0$ and expand $x\otimes y$ as the sum of four terms. It is enough to consider $x_0\otimes y_0$--term. In the $\os$--case, $x,y\in((\asa)^m)^-$ and $x_0, y_0$ may be assumed in $A^m_+$ and $B^m_+$, so that $x_0\otimes y_0\in C^m_+$. In the $\bU$--case, $x_0$ and $y_0$ are $\sigma$--weak limits from above of nets $(h_i), (k_j)$, where $h_i$ and $k_j$ are strongly lower semicontinuous. Then, since multiplication is separately $\sigma$--weakly continuous, to approximate $x_0\otimes y_0$ from above at a state $\varphi$, we first choose $j_0$ such that $\varphi(x_0\otimes(k_{j_0}-y_0))<{\ep\over 2} $; then choose $i_0$ such that $\varphi((h_{i_0}-x_0)\otimes k_{j_0})<{\ep\over 2}$. The rest is now clear. (One could also use Theorem 4 in the $\overline{S}$--, $\overline \bU_b$--, and $\bU_0$--cases.)

The slice maps were defined by J. Tomiyama \cite{T}. They take $a\otimes b$ into $f(a)b$ or $g(b)a$ for $f$ in $A^*$ or $g$ in $B^*$, and they extend to $\sigma$--weakly continuous maps from $C^{**}$ to $A^{**}$ or $B^{**}$. Since $g$ is a linear combination of positive functionals, it is enough to consider the case of a positive slice map. This plays the same role as the map $\psi$ in the proof of Lemma 2.
\enddemo

Of course, the slice maps factor through $(A\underset{\text{min}}\to\otimes B)^{**}$, so their usefulness is limited.

For $\bK$ the algebra of compact operators on $\ell^2$, $(A\otimes\bK)^{**}$ can be identified with $A^{**}\overline\otimes B(\ell^2)$, as is well known. The elements can be represented as infinite matrices over $A^{**}$, relative to the standard basis of $\ell^2$, and the only condition needed on such a matrix is that it represent a bounded operator.

\proclaim{Proposition 9}An element $x$ of $(A\otimes\bK)^{**}$ is in $\bU_\bC(A\otimes\bK)$ if and only if each matrix component of $x$ is in $\bU_\bC(A)$. The same is true for $\bB_s$, the sequential double--strong closure of the Borel algebra, and for the sequential double--strong closure of $\bU_{0\bC}$. Also, $x\in \bU_1(A\otimes\bK)$ if and only if each matrix component of $x$ is in $\bU_1(A)$. 
\endproclaim

\demo{Proof}Let $P_n$ in $B(\ell^2)$ be the projection on the span of the first $n$ basis vectors. If all components of $x$ are in $\bU_\bC(A)$, then $y_n=(1\otimes P_n)x(1\otimes P_n)$ is in $\bU_\bC(A\otimes\bK)$ by Proposition 8 (or by the proof of Proposition 3). Since $(y_n)$ converges double--strongly to $x$, and since $\bU_\bC$ is sequentially double--strongly closed by [P1, Theorem 3.7], we conclude $x\in \bU_\bC (A\otimes\bK)$. The converse follows from the last sentence of Proposition 8 (or from the proof of Proposition 3).

The second sentence is proved similarly.

For the last sentence we apply directly the definition of $\bU_1$ given in Remark 6. Note that the matrix multiplication formula, $(xy)_{ij}=\sum_k\, x_{ik}y_{kj}$, is a double--strongly convergent series, for $x, y$ in $(A\otimes\bK)^{**}$. Then it follows from the already proved $\bU$--case that if one of $x$ and $y$ is in $\bU_\bC(A\otimes\bK)$ and the other has all components in $\bU_1(A)$, then $xy\in \bU_\bC$. This shows one direction. For the other, take one of $x,y$ in $\bU_1(A\otimes\bK)$ and the other a diagonal matrix with components in $\bU_\bC(A)$.
\enddemo

We close with an open--ended problem. Note that Proposition 9 doesn't follow from Proposition 8, since some additional argument using the structure of $\bK$ was needed. Are there stronger results than Proposition 8 for tensor products? These could apply to fairly general tensor products (possibly one algebra should be nuclear and separable, for example) or to interesting special cases (assume one algebra is commutative, say).

\bye